\documentclass[10pt]{article}
\usepackage{a4,amsmath,amssymb,theorem,amsfonts,color}
\usepackage{multirow,graphicx,epsfig,float}
\usepackage{dsfont}

\title{Recursive low-rank approximation and model reduction of second-order systems}%

\author{Younes Chahlaoui%
        \thanks{Department of Mathematics,
                King Khalid University, KSA
                (\texttt{Younes.Chahlaoui@gmail.com})}}
\date{}

\mathcode`@="8000 
{\catcode`\@=\active\gdef@{\mkern1mu}}

\makeatletter
\def\mymatrix#1{\null\,\vcenter{\normalbaselines\m@th
    \ialign{\hfil$##$\hfil&&\quad\hfil$##$\hfil\crcr
      \mathstrut\crcr\noalign{\kern-\baselineskip}
      #1\crcr\mathstrut\crcr\noalign{\kern-\baselineskip}}}\,}
\makeatother
\def\mybmatrix#1{\left[ \mymatrix{#1} \right]}

\newcommand{\beq}{\begin{equation}}
\newcommand{\eeq}{\end{equation}}
\newcommand{\beqno}{\begin{equation*}}
\newcommand{\eeqno}{\end{equation*}}
\newcommand{\eqn}{\begin{eqnarray}}
\newcommand{\eeqn}{\end{eqnarray}}
\newcommand{\eqnn}{\begin{eqnarray*}}
\newcommand{\eeqnn}{\end{eqnarray*}}
\newcommand{\bma}{\left[\begin{array}}
\newcommand{\ema}{\end{array}\right]}
\newcommand{\bsy}{\left \{ \begin{array}}
\newcommand{\esy}{\end{array}\right.}
\newcommand{\dfr}[2]{\displaystyle\frac{#1}{#2}}

\def\norm2#1{\|{#1}\|_2}

\def\H2{\mathcal{H}_2}
\def\Hinf{\mathcal{H}_\infty}

\def\R{\mathbb{R}}

\let\oldref\ref
\def\ref#1{{\normalfont\oldref{#1}}}
\def\eqref#1{{\normalfont(\oldref{#1})}}
{\theorembodyfont{\it}
\newtheorem{theorem}{Theorem}[section]}
{\theorembodyfont{\it}
}
{\theorembodyfont{\it}
}
{\theorembodyfont{\it}
}
 {\theorembodyfont{\it}
 }
{\theorembodyfont{\it}
}

\newtheorem{algorithm}[theorem]{Algorithm}

\newcounter{mylineno}
\makeatletter
\let\oldtabcr\@tabcr

\def\mynewline{\refstepcounter{mylineno}%
               \llap{\footnotesize\arabic{mylineno}\hspace{5pt}}%
              }

\gdef\@tabcr{\@stopline \@ifstar{\penalty%
           \@M \@xtabcr}\@xtabcr\mynewline}

\newenvironment{code}{%
                        \mathcode`\:="603A  
                        
                        \setcounter{mylineno}{0}
                        \par
                        \upshape
                        \begin{list} 
                           {} {\leftmargin = 1cm}
                        \item[]
                        \begin{tabbing}

                           \hspace*{.3in} \= \hspace*{.3in} \=
                           \hspace*{.3in} \= \hspace*{.3in} \= \kill
                           \mynewline
                       }{\end{tabbing}\end{list}}
\makeatother

\begin{document}
\maketitle
\begin{abstract}
We present an adaptation of two recent low-rank approximation technique proposed for first-order model reduction systems to the second-order systems.
The resulting reduced order models are guaranteed to keep the second order structure which has a physical meaning. The quality of the approximation is shown using numerical simulations on benchmark examples.
\end{abstract}


\section{Introduction}\label{INTRO}

Model reduction plays an important role in simulation and control of large dynamical systems \cite{Ant05,BMS05,Schi2008}. For second-order dynamical systems
\begin{equation}\label{eq1}
M\ddot{q}(t) + D\dot{q}(t) + Kq(t) = Fu(t),
\quad y(t)=Gq(t),
\end{equation}
with inputs $u(t)\in\R^m$, states $q(t)\in\R^N$ and outputs $y(t)\in \R^p$ (all matrices of appropriate dimensions), the problem is more challenging as these systems have a systemic constraint. Indeed, these systems arise naturally in many engineering applications in general nonlinear and physical (see \cite{Gaw08} and the references in).

To handle these systems, we often linearize and reformulate the original system by introducing $x(t)^T=\begin{bmatrix}q(t)^T & \dot{q}(t)^T\end{bmatrix}$, which will lead to a first-order formulation (also called state-space description) of \eqref{eq1} as
\begin{equation}\label{eq2}
E\dot{x}(t) =Ax(t)+Bu(t),
\quad y(t)=CX(t),
\end{equation}
where $E = \mybmatrix{W&0\cr 0& M}, \quad A = \mybmatrix{0&W\cr -K&-D}, \quad B = \mybmatrix{0\cr F}, \quad C = \mybmatrix{G&0}$ and $W$ any nonsingular matrix ($W=I$ in practice). The structure in the quadruplet $\{E, A, B, C\}$ is significantly associated to the physical origin of the variables $q(t)$ in \eqref{eq1}. Consequently, in order to understand what the reduced model will produce as results, we need to consider mainly (second-order) structure preserving model reduction methods.

Various techniques have been developed in the literature for model reduction of second-order systems. We can classify these techniques (vaguely) in two families: the engineers's methods and the control theorists' methods. The first class considers the second order form directly and based on some engineering understanding (often the notion of master/slave variables is used) is tailoring a reduction to each specific case. This case by case approach prohibit most system analyses. Moreover, very often the reduction is based on the choice of some particular characteristic frequencies lying for example in a particular range. A critical drawback of this approach is that we only choose half of these frequencies as the second half is automatically obtained (complex conjugate of the first half). These frequencies not chosen are often sources of hazardous behavior.  Guyan reduction \cite{Gaw08,Guy65} and Component Mode Synthesis (CMS) \cite{CB68,MS02} belong to this family of methods and are frequently used in commercial FEM software packages. Other techniques are dynamic reduction \cite{SV07}, Improved Reduction System (IRS) \cite{FGP95,MS02} and System Equivalent Expansion Reduction Process (SEREP) \cite{OAR89}.

The control theorists' methods consist of second order structure preserving model reduction techniques. Within this family, we can distinguish between three classes. The first class deals with the second order directly without the need of converting it into a state space representation. The second order Krylov subspace methods are an example of such direct second order methods \cite{BH03,BMS05,Schi2008,MS96,SL07,Ten04,Van04}. Another class is the indirect methods represented by the second order balanced truncation methods \cite{CLMVPVD02,CLVPVD04b,CLVPVD04a,MS96,RS08}. In these methods, second order Gramians are defined by making a variable correspondences between second order form and first-order form. All above methods still have many drawbacks. For
instance, neither a global error bound is guaranteed to be obtained, nor the stability is preserved in the reduced systems. Moreover, as a sufficient condition, only block diagonally projectors have to be chosen in order to satisfy the structural conditions. The third class is based on some remarks made first in \cite{MS96} (up to the knowledge of the authors).
The key idea is to come up with a reduced system, via a first-order method, and be able to unfold it into a second order form. In \cite{MS96}, different algebraic conditions, not necessarily easy to check, were presented. A sufficient and necessary condition was presented in \cite{CH2015}. It is based on a characteristic property for all second order systems via their first Markov parameter. This property guarantees that the first-order form is a linearization of a second order form.

In this note we present  an adaptation of some recent low-rank approximation techniques proposed for the first-order case \cite{Cha2010} to the second order formulation, in such a way that the resulting reduced order model is in the second order formulation directly.

The note is as follows below. Section 2 is dedicated to differential and differences second-order systems. In Section 3 we present the modelling and projection of dynamics of second-order systems. Both recursive low-rank methods (RLRG and RLRH) are presented in Section 4. Section 5 is dedicated to numerical examples and a comparison of our methods with some other methods for model reduction of first-order and second-order systems. We conclude with some comments and open problems.

\section{Differential and differences second-order systems}
The RLRG and RLRH algorithms are designed for discrete systems. We need first to transform the differential equations in \eqref{eq1} using differences. For this, we can approximate the first and second differentials using Newton difference formula for example. Table \ref{tab1} is showing different relations for doing these transformations, where $h$ is the step size and $q(ih)=q_i$.

\begin{table}
\begin{center}
\begin{tabular}{|c|c|l|}
  \hline
  $\ddot{q}(t)$ & $\dot{q}(t)$ &  \\\hline
   & $\dfr{q_{i+1}-q_i}{h}$ &$\begin{array}{lcl}\bar{M}&=&(M+hD)/h^2\\ \bar{D}&=&(h^2K-2M-hD)/h^2\\ \bar{K}&=&M/h^2\end{array}$\\\cline{2-3}
  $\dfr{q_{i+1}-2q_i+q_{i-1}}{h^2}$ & $\dfr{q_i-q_{i-1}}{h}$ &  $\begin{array}{lcl}\bar{M}&=&M/h^2\\ \bar{D}&=&(h^2K-2M+hD)/h^2\\ \bar{K}&=&(M-hD)/h^2\end{array}$\\\cline{2-3}
   & $\dfr{q_{i+1}-q_{i-1}}{2h}$ &  $\begin{array}{lcl}\bar{M}&=&(2M+hD)/2h^2\\ \bar{D}&=&(h^2K-2M)/h^2\\ \bar{K}&=&(2M-hD)/2h^2\end{array}$\\
  \hline
\end{tabular}
\caption{Relation between differential and difference second-order systems}\label{tab1}
\end{center}
\end{table}

Remark that if $M$, $D$, and $K$ are symmetric, $\bar{M}$, $\bar{D}$, and $\bar{K}$ will also be symmetric. We obtain an equivalent difference equation representing also a second-order system
 \begin{equation}\label{eq3}
\bar{M}q_{i+1}+\bar{D}q_i+\bar{K}q_{i-1} = Fu_i,
\quad y_i=Gq_i.
\end{equation}
The transfer functions associated with \eqref{eq1} and \eqref{eq3} are
\[T_f(s)=GP(s)^{-1}F,\quad\textrm{and}\quad T_f(z)=GP(z)^{-1}F,\]
where $P(s)=Ms^2+Ds+K$, and $P(z)=\bar{M}z+\bar{D}+\bar{K}z^{-1}$ are the characteristic polynomial matrices. The zeros of $\det(P(\cdot))$ are also known as the characteristic frequencies of the system and play an important role in engineering model reduction techniques. Stability of the system, e.g., implies that these zeros must lie in the stability region, which is the left half plane for the continuous-time case and the open unit desk of the complex plane (the open disk centered at 0 of radius 1) for the discrete-time case.
In the sequel, for simplification we will use the quintuplet $\{M, D, K, F, G\}$ for both continuous and discrete case.

\section{Modelling and projection of dynamics of second-order systems}

Most interesting model reduction methods are those based on the projection of dynamics. We seek to construct two projection matrices $X$ and $Y\in \R^{2N\times 2n}$ ($n\ll N$) ($Y^TX=I_{2n}$ and $XY^T$ is a projector). We consider the following partition, where each block is $N\times n$
\[X=\mybmatrix{X_{11}&X_{12}\cr X_{21}&X_{22}},\quad Y=\mybmatrix{Y_{11}&Y_{12}\cr Y_{21}&Y_{22}}.\]
Starting from any linearization \eqref{eq2}, to preserve the special structure \eqref{eq2} we have to choose $X$ and $Y$ so that
\[Y^T\mybmatrix{I&0\cr 0& M}X=\mybmatrix{T_1&0\cr 0& \hat{M}}, \quad Y^T\mybmatrix{0&I\cr -K&-D}X=\mybmatrix{0&T_2\cr -\hat{K}&-\hat{D}},\]
 \[Y^T\mybmatrix{0\cr F}=\mybmatrix{0\cr \hat{F}}\quad\textrm{and}\quad \mybmatrix{G&0}X=\mybmatrix{\hat{G}&0}\]
Where $T_i$, $i=1,2$ are nonsingular matrices.

Sufficient conditions to obtain this for all $M$, $D$, $K$, $F$ and $G$, are to choose $X$ and $Y$ block diagonal, i.e.
\[X_{12}=X_{21}=Y_{12}=Y_{21}=0\]
provided $T_1=Y_{11}^TX_{11}$ and $T_2=Y_{11}^TX_{22}$ are nonsingular \cite{CLVPVD04b}.
In such a case,
\[\hat{M}=Y_{22}^TMX_{22}, \hat{D}=Y_{22}^TDX_{22}, \hat{K}=Y_{22}^TKX_{11}, \hat{F}=Y_{22}^TF, \hat{G}=GX_{11}.\]
We may choose $\tilde{X}=XT^{-1}$ where $T=blkdiag(T_1, T_2)$ to obtain a reduced order model in standardized form (i.e., $\hat{E}=I$).
The reduced order model matrices are now
\[\hat{M}=Y_{22}^TM\tilde{X}_{22}, \hat{D}=Y_{22}^TD\tilde{X}_{22}, \hat{K}=Y_{22}^TK\tilde{X}_{11}, \hat{F}=Y_{22}^TF, \hat{G}=G\tilde{X}_{11}.\]
In the sequel we show how to compute such projection matrices using low-rank approximations of the Gramians of first-order formulation.

\section{RLRG and RLRH for second-order systems}

The Recursive Low-Rank Gramian (RLRG) and Recursive Low-Rank Hankel (RLRH) algorithms presented in \cite{Cha2010} can be adapted easily for second-order systems (at least for discretized systems). The idea is to use the special structure of the matrices $A$, $B$ and $C$ of these systems to obtain a recurrence involving directly the matrices of the second-order system $M$, $D$, $K$, $F$ and $G$. Since we have a second-order system, one may expect that the procedure will involve a recurrence on two time steps on a matrix representing the dominant subspace. From this matrix, we will construct the projection matrices directly for the second-order system.

It is pointed out that the updated version of the dominant subspace for both algorithms (RLRG and RLRH) is obtained via the recursion
\eqn S_n(i+1)&=&\underbrace{\bma{c|c} AS_n(i) & B\ema}_{M_1} \bma{c}V_{c_1}\\\hline V_{c_2}\ema, \textrm{and} \\ R_n(i+1)&=&\underbrace{\bma{c|c} A^TR_n(i) & C^T\ema}_{M_2} \bma{c}V_{o_1}\\\hline V_{o_2}\ema,\label{eq4}\eeqn
Where $V_{c_1}, V_{o_1} \in \mathbb{R}^{n\times n}$ come from the {SVD} of $M_1$ and $M_2$, respectively for the RLRG algorithm or from the SVD of $M_2^TM_1$ for the RLRH algorithm \cite{Cha2010}.

For a second-order system, one may consider that $S_\bullet(\cdot)$ corresponds in fact a collection of two consecutive versions of the dominant subspace of the second-order system, i.e., for example for $S_n(\cdot)$
\[S_n(i+1)=\bma{c}S_n^s(i)\\ S_n^s(i+1)\ema,\]
(subscript $s$ referring to ``second-order system"). It follows from this and from \eqref{eq4} that we have
\beq \bma{c}S_n^s(i)\\S_n^s(i+1)\ema=\bma{c|c}\bma{cc} 0&I\\ -M^{-1}K& -M^{-1}D\ema\bma{c}S_n^s(i-1)\\S_n^s(i)\ema& \bma{c}0\\ M^{-1}F\ema\ema \bma{c}V_1\\\hline V_2\ema.\label{eq5}\eeq
Note first that $S_n^s(i)$ on the left hand side is not the same as on the right hand side, but this is an intermediate updating version which may allow us to update the whole subspaces for state-space description. If we denote it $S_{n_+}^s(i)$, \eqref{eq5} can be rewritten as
\beq\bsy{ccl} S_{n_+}^s(i)&=&S_n^s(i)V_1,\\ S_{n}^s(i+1)&=&-M^{-1}KS_n^s(i-1)V_1-M^{-1}DS_n^s(i)V_1+M^{-1}FV_2.\esy\label{eq6}\eeq
Here, the first equation is redundant and one may consider just the second one. But, actually this is not correct as it is an update of $S_n(i)$ which will affect not $S_n(i+1)$ but $S_n(i+2)$ in the next step.

For $R_n(\cdot)$, the adaptation is more complicated. We have
\beq \bma{c}R_n^s(i)\\R_n^s(i+1)\ema=\bma{c|c}\bma{cc} 0&-K^TM^{-T}\\ I& -D^TM^{-T}\ema\bma{c}R_n^s(i-1)\\R_n^s(i)\ema& \bma{c}0\\ G^T\ema\ema \bma{c}U_1\\\hline U_2\ema.\label{eq7}\eeq
which yields
\beq\bsy{ccl} R_{n_+}^s(i)&=&-K^TM^{-T}R_n^s(i)U_1,\\ R_{n}^s(i+1)&=&R_n^s(i-1)U_1-D^TM^{-T}R_n^s(i)U_1+G^TU_2.\esy\label{eq8}\eeq
Here, the first equation is not redundant and we have to take it into consideration. One can remark also that the first equation will affect the dominant subspace only after two steps.

Using recurrences \eqref{eq6} and \eqref{eq8}, we obtain after a certain number of iterations, two $N\times n$ matrices $S_n$ and $R_n$. These matrices are used to construct the projection matrices $X$ and $Y$ (previous section) as we will show later. For the convergence and accuracy of both algorithms RLRG and RLRH we refer to \cite{Cha2010} for a more detailed analysis.

In the sequel we present the two algorithms adapted for second-order systems, namely Second-order Recursive Low-Rank Gramian (SRLRG) and Second-order Recursive Low-Rank Hankel (SRLRH).
To simplify the algorithms let us define the following matrices
\[M_1(i)=\bma{c|c}\bma{cc}0&I\\ -M^{-1}K&-M^{-1}D\ema\bma{c}S_n(i-1)\\ S_n(i)\ema&\bma{c}0\\ M^{-1}F\ema\ema,\]
\[M_2(i)=\bma{c|c}\bma{cc}0&-K^TM^{-T}\\ I&-D^TM^{-T}\ema\bma{c}R_n(i-1)\\ R_n(i)\ema&\bma{c}0\\ G^T\ema\ema,\]

\begin{algorithm}\label{SRLRG}
Second-order Recursive Low-Rank Gramians algorithm (SRLRG):\\
Given $\{M, D, K, F, G\}$ and $n$, we compute $\{\hat{M}, \hat{D}, \hat{K}, \hat{F}, \hat{G}\}$
\end{algorithm}
\begin{code}
  Let $S_n(0)$, $S_n(1)$, $R_n(0)$, $R_n(1)$ be any $N\times n$ random matrices\\
  for $i = 1:\tau$\\
  \> Compute the SVDs $M_1(i)=U_c\Sigma_c V_c^T$, and $M_2(i)=U_o\Sigma_o V_o^T$,\\
  \> Let $V_c(:,1:n)=\bma{c}V_{c_1}\\ V_{c_2}\ema$, and $V_o(:,1:n)=\bma{c}V_{o_1}\\ V_{o_2}\ema$, where $V_{c_1}, V_{o_1}\in \R^{n\times n}$\\
  \> Compute the new $S_n(\cdot)$ and $R_n(\cdot)$\\
    \> \> $S_{n_+}=S_n(i)V_{c_1},$\\
    \> \> $S_n(i+1)=-M^{-1}KS_n(i-1)V_{c_1}-M^{-1}DS_n(i)V_{c_1}+M^{-1}FV_{c_2}$\\
     \> \> and\\
    \> \> $R_{n_+}=-K^TM^{-T}R_n(i)V_{o_1},$\\
    \> \> $R_n(i+1)=R_n(i-1)V_{o_1}-D^TM^{-T}R_n(i)V_{o_1}+G^TV_{o_2}$\\
  end\\
  Compute the SVD $S_n(\tau)^TR_n(\tau)=U\Sigma V^T$\\
  Construct $X=S_n(\tau)U\Sigma^{-1/2}$, and $Y=R_n(\tau)V\Sigma^{-1/2}$\\
  Construct $\hat{M}=Y^TMX$, $\hat{D}=Y^TDX$, $\hat{K}=Y^TKX$, $\hat{F}=Y^TF$, $\hat{G}=GX$.
\end{code}

\begin{algorithm}\label{SRLRH}
Second-order Recursive Low-Rank Hankel algorithm (SRLRH):\\
Given $\{M, D, K, F, G\}$ and $n$, we compute $\{\hat{M}, \hat{D}, \hat{K}, \hat{F}, \hat{G}\}$
\end{algorithm}
\begin{code}
  Let $S_n(0)$, $S_n(1)$, $R_n(0)$, $R_n(1)$ be any $N\times n$ random matrices\\
  for $i = 1:\tau$\\
  \> Compute the SVD:\quad $M_2^TM_1=U\Sigma V^T$,\\
  \> We use the following decomposition (where $V_1$, $U_1 \in \R^{n\times n}$)\\ \> \>$V(:,1:n)=\bma{c}V_1\\ V_2\ema$, and $U(:,1:n)=\bma{c}U_1\\ U_2\ema$ \\
  \> Compute the new $S_n(\cdot)$ and $R_n(\cdot)$:\\
    \> \> $S_{n_+}=S_n(i)V_1,$\\
    \> \> $S_n(i+1)=-M^{-1}KS_n(i-1)V_1-M^{-1}DS_n(i)V_1+M^{-1}FV_2$\\
     \> \> and\\
    \> \> $R_{n_+}=-K^TM^{-T}R_n(i)U_1,$\\
    \> \> $R_n(i+1)=R_n(i-1)U_1-D^TM^{-T}R_n(i)U_1+G^TU_2$\\
  end\\
  Compute the SVD:\quad $S_n(\tau)^TR_n(\tau)=U\Sigma V^T$\\
  Construct\quad $X=S_n(\tau)U\Sigma^{-1/2}$, and $Y=R_n(\tau)V\Sigma^{-1/2}$\\
  Construct\quad $\hat{M}=Y^TMX$, $\hat{D}=Y^TDX$, $\hat{K}=Y^TKX$, $\hat{F}=Y^TF$, $\hat{G}=GX$.
\end{code}

A couple of comments are in order here. First, the choice of the reduced order $n$ can be chosen a priori by the user or chosen by the method following specification of a tolerance by the user who still have to specify a bounding interval for it. Second the for loop used in both algorithms needs the specification of $\tau$ the total number of iterations. It was shown in \cite{Cha2010} that in the worst case (when the system is close to instability) one needs iterations three times the dimension of the original system to get a kind of invariance upon the principal subspaces computed by the SVDs. We can replace the for loop by a while loop based on a stopping criterion build on the invariance upon the principal subspaces computed by the SVDs. In this case, one has only to specify a tolerance value on the angle between old and new principal subspaces (see \cite[Section 7]{Cha2010}).

\section{Numerical Examples}\label{examples}

This section makes a numerical comparison of several procedures for model order reduction for second-order systems. The considered test methods are:
\begin{itemize}
\item the classic Balanced Truncation (\textbf{BT}): one of the best model reduction methods, used here as a quality's reference,
 \item Mayer and Srinivasan methods: Free Velocity method (\textbf{FV}), Zero Velocity method (\textbf{ZV}) and (\textbf{MS}) for the case where we apply free and zero velocity in the same time \cite{MS96},
 \item Su and Craig (\textbf{SC}) method: one of the best Krylov method,
     \item the Second-order Balanced Truncation (\textbf{SoBT}) \cite{CLVPVD04b},
     \item our methods: Second-order Recursive Low-Rank Gramian (\textbf{SoRLRG}) and Second-order Recursive Low-Rank Hankel (\textbf{SoRLRH}).
 \end{itemize}
These algorithms are applied to benchmark models from \cite{CPVD05}: the Building model, the CDplayer model, and ISS 1R and 12A models. Table. \ref{tab-data} is showing the data about each model.

\begin{table}[H]
\begin{center}
\begin{tabular}[H]{|c|c|c|c|c|}
\hline &$2N$&$m$&$p$&$2n$\\
\hline build model&48&1&1&10\\
\hline ISS 1R model&270&3&3&32\\
\hline ISS 12A model&1412&3&3&196\\
\hline
\end{tabular}
\caption{Summary of data of the benchmark models.}\label{tab-data}
\end{center}
\end{table}

\begin{center}
\begin{figure}[H]
\begin{center}
\begin{minipage}[H]{.9\linewidth}
\centering{\epsfig{figure=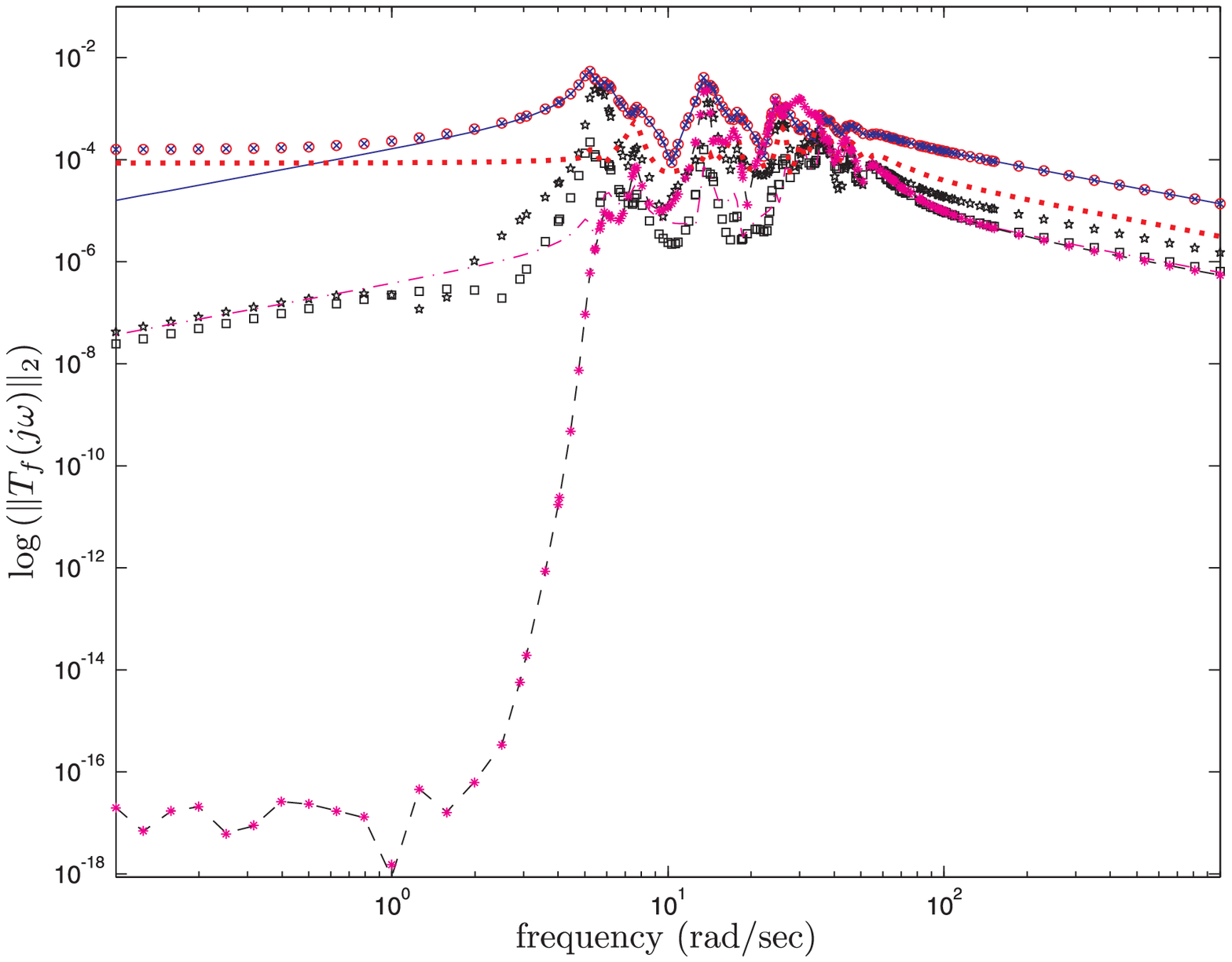,width=.9\linewidth,height=7cm}}
\caption[$\sigma_{\max}$-plot of the frequency responses for
Building model.]{\centering{$\sigma_{\max}$-plot of the frequency
responses for Building model.}}\label{build7}
\end{minipage}

\medskip

\begin{minipage}[H]{.9\linewidth}
\centering{\epsfig{figure=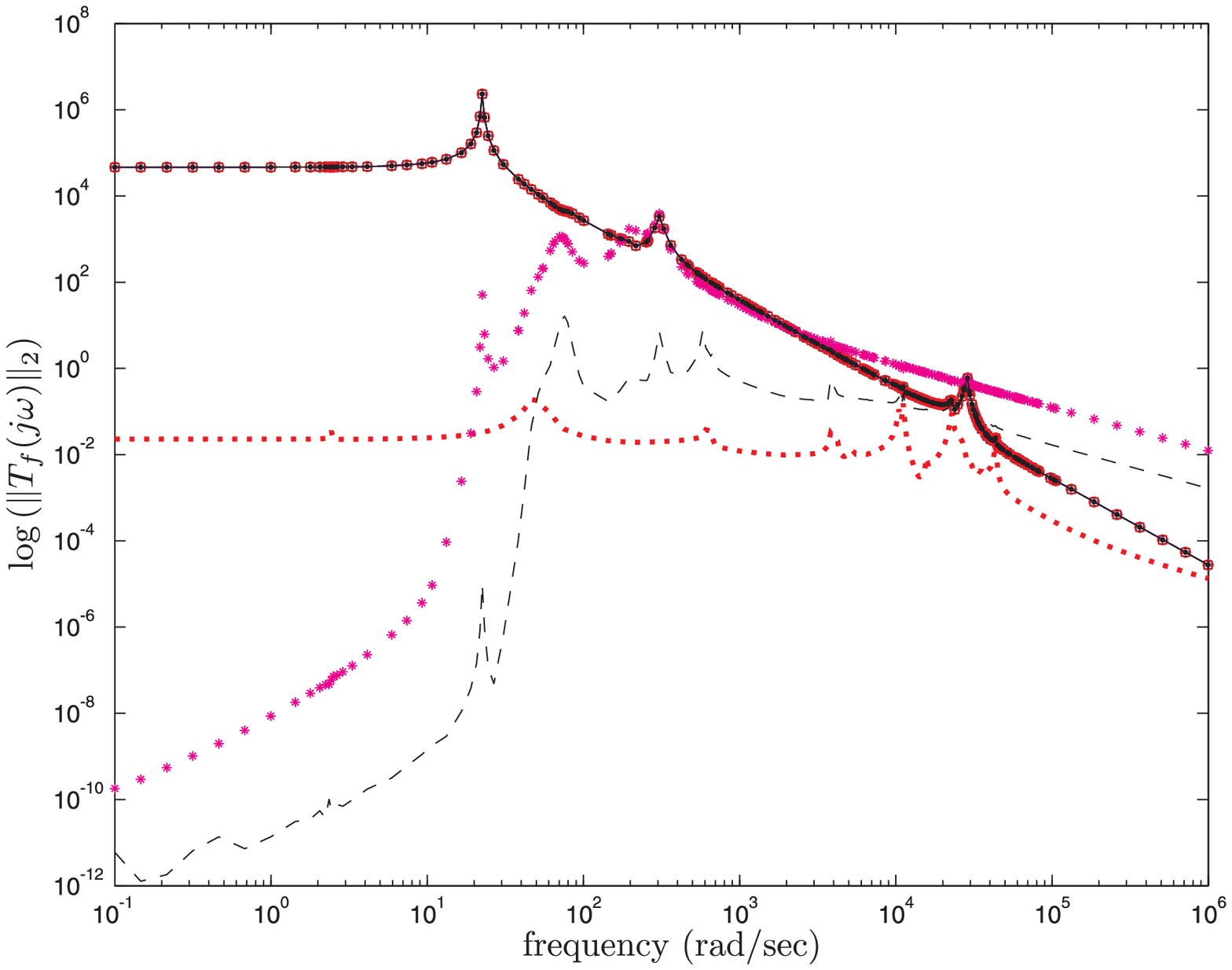,width=.9\linewidth,height=7cm}}
\caption[$\sigma_{\max}$-plot of the frequency responses for
CDplayer model.]{\centering{$\sigma_{\max}$-plot of the frequency
responses for CDplayer model.}}\label{cd7}
\end{minipage}

\medskip

\textsc{Legend~:} \textcolor{blue}{\rule[.7mm]{3mm}{.2mm}} full
model, \textcolor{red}{$\cdots$} BT error system, \textcolor{black}{$\square$} SOBT error system,\\
\textcolor{red}{o} Mayer error system, \textcolor{blue}{$\times$}
Su error system, \textcolor{black}{$\star$} FV error system,
\textcolor{black}{\rule[.7mm]{1mm}{.2mm}\!$\cdot$\!\rule[.7mm]{1mm}{.2mm}}
ZV error system,\\ \rule[.7mm]{1mm}{.2mm}\;\rule[.7mm]{1mm}{.2mm}
SRLRG error system, \textcolor{magenta}{$\ast\ast\ast$} SRLRH
error system.
\end{center}
\end{figure}
\end{center}

\begin{center}
\begin{figure}[H]
\begin{center}
\begin{minipage}[H]{.9\linewidth}
\centering{\epsfig{figure=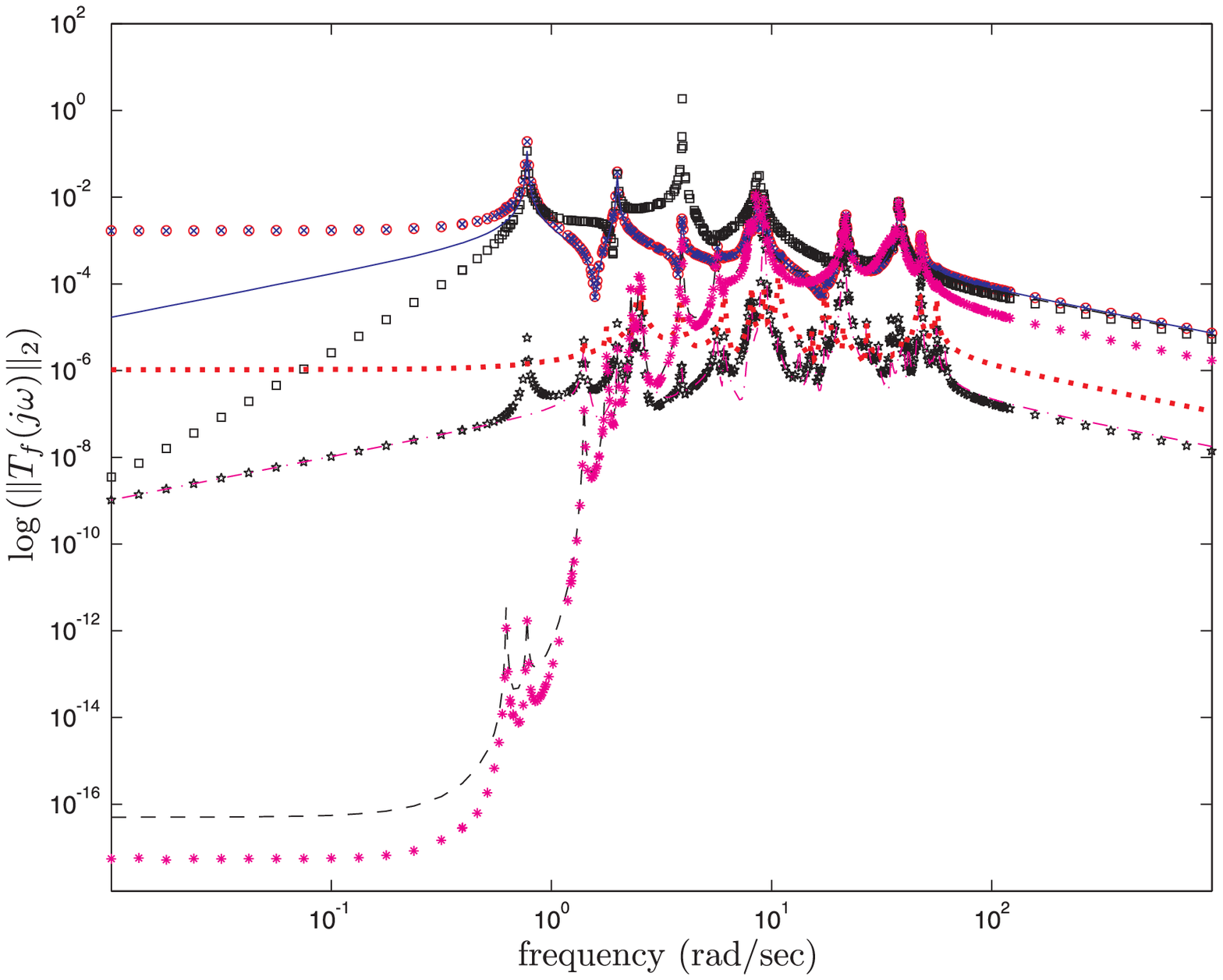,width=.9\linewidth,height=7cm}}
\caption[$\sigma_{\max}$-plot of the frequency responses for ISS
1R model.]{\centering{$\sigma_{\max}$-plot of the frequency
responses for ISS 1R model.}}\label{iss17}
\end{minipage}

\medskip

\begin{minipage}[H]{.9\linewidth}
\centering{\epsfig{figure=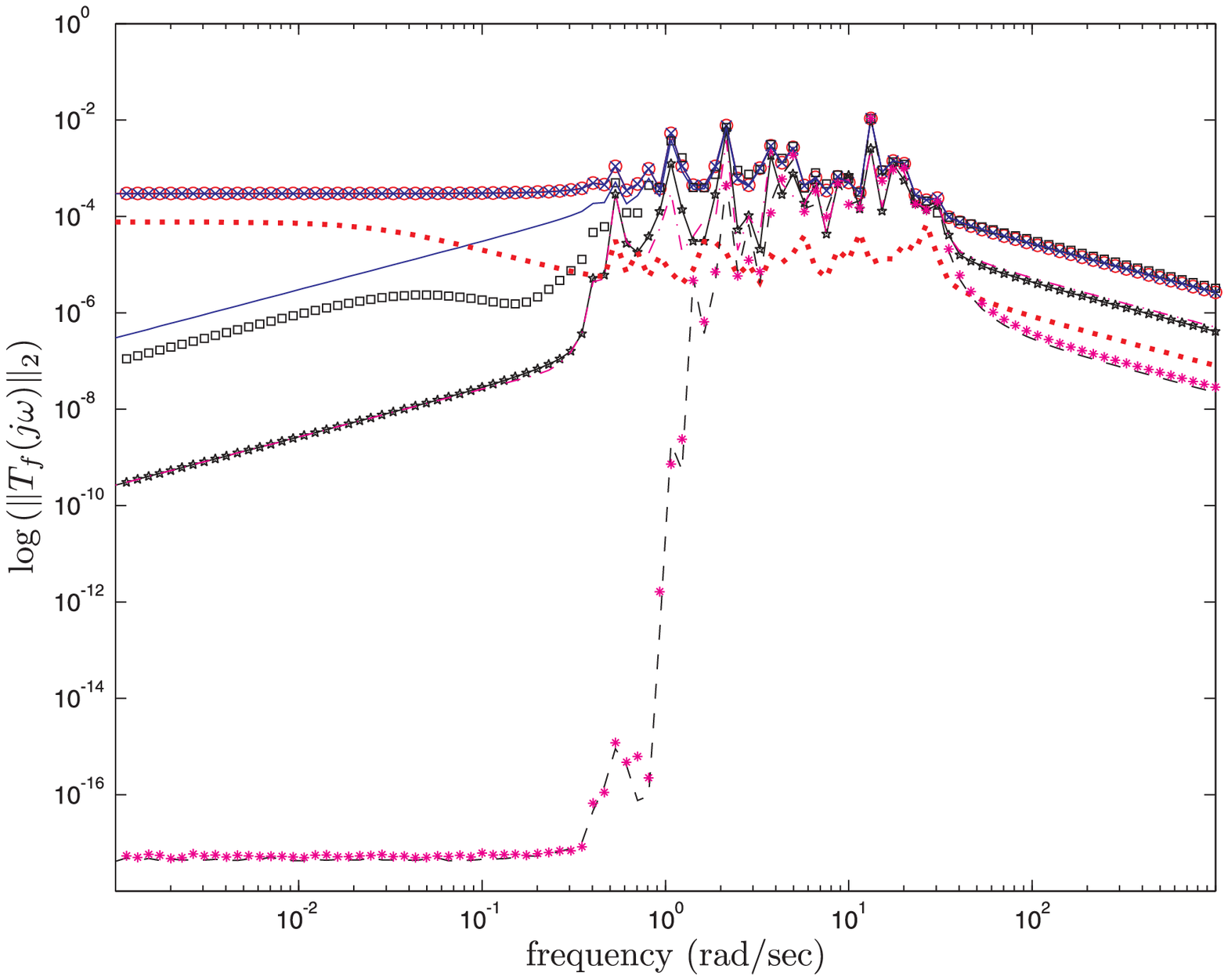,width=.9\linewidth,height=7cm}}
\caption[$\sigma_{\max}$-plot of the frequency responses for ISS
12A model.]{\centering{$\sigma_{\max}$-plot of the frequency
responses for ISS 12A model.}}\label{iss27}
\end{minipage}

\medskip

\textsc{Legend~:} \textcolor{blue}{\rule[.7mm]{3mm}{.2mm}} full
model, \textcolor{red}{$\cdots$} BT error system, \textcolor{black}{$\square$} SOBT error system,\\
\textcolor{red}{o} Mayer error system, \textcolor{blue}{$\times$}
Su error system, \textcolor{black}{$\star$} FV error system,
\textcolor{black}{\rule[.7mm]{1mm}{.2mm}\!$\cdot$\!\rule[.7mm]{1mm}{.2mm}}
ZV error system,\\ \rule[.7mm]{1mm}{.2mm}\;\rule[.7mm]{1mm}{.2mm}
SRLRG error system, \textcolor{magenta}{$\ast\ast\ast$} SRLRH
error system.
\end{center}
\end{figure}
\end{center}

The comparison is made on the basis of the relative reduction error measured according to the $\Hinf$ norm, i.e., the ratio of the $\Hinf$ norms of the "error" system $\mathcal{S}-\mathcal{S}_n$ and the full system $\mathcal{S}$, i.e.,
\[\textbf{rre}(method)=\dfr{\|\mathcal{S}-\mathcal{S}_n^{method}\|_\infty}{\|\mathcal{S}\|_\infty}.\]
The results are shown in Table. \ref{Hinfresults}. In Figures~\ref{build7}, \ref{cd7}, \ref{iss17}, and \ref{iss27}
we show the $\sigma_{\max}$-plot of the frequencies of the full
system and the error system for each benchmark model.

The method of Su and Craig produces for the CD player model an unstable reduced-order model, but in general we have remarked that even if the original system is strictly stable, the reduced-order model obtained via this method moves to a region close to instability. The methods of Meyer seem to work quite well, but in general the method SOBT behaves better. The adaptation of RLRG and RLRH, namely SoRLRG and SoRLRH yield the best results. Note
that we still use a stopping criterion based on the tolerance values given in \cite{Cha2010}.

\begin{table}[H]
\begin{center}
\begin{tabular}[H]{|c|c|c|c|c|}
\hline&&&&\\
&Building&CDplayer&ISS 1R&ISS 12A\\&model&model&model& model\\
\hline
$\|\mathcal{S}\|_\infty$&0.0053&$2.3198.10^6$&0.1159&0.0107\\
\hline $rre(BT)$&0.1143&$8.0704.10^{-8}$&0.0013&0.0071\\
\hline $rre(MS)$&1.0181&1.0000&1.6324&1.0028\\
\hline $rre(SC)$&1.0182&unstable&1.6324&1.0063\\
\hline $rre(SoBT)$&0.0424&1.0000&0.1054&0.9966\\
\hline $rre(FV)$&0.4509&1.0000&0.0548&0.5659\\
\hline $rre(ZV)$&0.0379&1.0000&0.0513&0.3537\\
\hline $rre(SoRLRG)$&0.4301&$6.8931.10^{-6}$&0.1023&0.9697\\
\hline $rre(SoRLRH)$&0.4320&$1.7.10^{-6}$&0.0979&0.9390\\
\hline
\end{tabular}
\caption{$\mathcal{H}_\infty$ norm of benchmark models, and the
error systems.}\label{Hinfresults}
\end{center}
\end{table}

\section{Concluding remarks}\label{concrem}

In this paper, we presented an adaptation of two recent low-rank approximation technique proposed for first-order model reduction systems to the second-order systems. The resulting reduced order models are guaranteed to keep the second order structure which has a physical meaning. The quality of the approximation is shown using numerical simulations on benchmark examples. Unlike all other methods in the literature, the reduced order models produced by our methods are guaranteed to be stable.

Despite the obviously desirable features of the approach proposed here, many open questions remain. There are a number of refinements with respect to performance, convergence, and accuracy which require more theoretical and algorithmic analysis. There are two particularly interesting features concerning the stability of the obtained reduced system and bounds on the quality.

\section*{Acknowledgements}

I gratefully acknowledge the helpful remarks and suggestions of Paul Van Dooren and Danny Sorensen.

\bibliographystyle{acm}
\bibliography{SecOrd1}

\begin{thebibliography}{10}

\bibitem{Ant05}
{\sc Antoulas, A.~C.}
\newblock {\em Approximation of Large-Scale Dynamical Systems}.
\newblock Society for Industrial and Applied Mathematics, Philadelphia, PA,
  USA, 2005.

\bibitem{BH03}
{\sc Bastian, J., and Haase, J.}
\newblock Order reduction for second order systems.
\newblock In {\em Proceedings 4th MATHMOD Vienna\/} (2003), I.~Troch and
  F.~Breitenecker, Eds., pp.~418--424.

\bibitem{BMS05}
{\sc Benner, P., Mehrmann, V., and Sorensen, D.~C.}, Eds.
\newblock {\em {Dimension reduction of large-scale systems}}, vol.~45 of {\em
  Lecture Notes in Computational Science and Engineering}.
\newblock Spring{\-}er-Ver{\-}lag, Berlin, 2005.

\bibitem{Cha2010}
{\sc Chahlaoui, Y.}
\newblock Two efficient {SVD/Krylov} algorithms for model order reduction of
  large scale systems.
\newblock {\em Electronic Transactions on Numerical Analysis (ETNA) 38\/}
  (2011), 113--145.
\newblock MIMS Eprint 2010.11.

\bibitem{CH2015}
{\sc Chahlaoui, Y., and Al-Qasmi, H.}
\newblock A new structure preserving model reduction method for second-order
  systems.
\newblock {\em submitted to SIAM J. Control and Optimization\/}.

\bibitem{CLMVPVD02}
{\sc Chahlaoui, Y., Lemonnier, D., Meerbergen, K., Vandendorpe, A., and {Van
  Dooren}, P.}
\newblock Model reduction of second order systems.
\newblock Paper 26984-4, MTNS2002, Notre Dame, USA, 2002.

\bibitem{CLVPVD04b}
{\sc Chahlaoui, Y., Lemonnier, D., Vandendorpe, A., and {Van Dooren}, P.}
\newblock Second order balanced truncation.
\newblock {\em Linear Algebra Appl. 415}, 2-3 (2004), 373--384.

\bibitem{CLVPVD04a}
{\sc Chahlaoui, Y., Lemonnier, D., Vandendorpe, A., and {Van Dooren}, P.}
\newblock Second order structure preserving balanced truncation.
\newblock Paper TA8-6, MTNS04, Leuven, 2004.

\bibitem{CPVD05}
{\sc Chahlaoui, Y., and {Van Dooren}, P.}
\newblock {\em Benchmark examples for model reduction of linear time invariant
  dynamical systems}.
\newblock Vol.~45 of Benner et~al. \cite{BMS05}, 2005, pp.~379--392.

\bibitem{CB68}
{\sc Craig, R.~R., and Bampton, M. C.~C.}
\newblock Coupling of substructures for dynamic analyses.
\newblock {\em AIAA Journal 6\/} (1968), 1313--1319.

\bibitem{FGP95}
{\sc Friswell, M.~I., Garvey, S.~D., and Penny, J. E.~T.}
\newblock Model reduction using dynamic and iterated {IRS} techniques.
\newblock {\em J. Sound Vib. 186\/} (1995), 311--323.

\bibitem{Gaw08}
{\sc Gawronski, W.}
\newblock {\em Modeling and control of antennas and telescopes.}
\newblock Mechanical Engineering Series. Spring{\-}er-Ver{\-}lag, New York,
  2008.

\bibitem{Guy65}
{\sc Guyan, J.}
\newblock Reduction of stiffness and mass matrices.
\newblock {\em AIAA Journal 3}, 2 (1965), 380.

\bibitem{MS96}
{\sc Meyer, D.~G., and Srinivasan, S.}
\newblock Balancing and model reduction for second-order form linear systems.
\newblock {\em IEEE Trans. Autom. Control. 41}, 11 (1996), 1632--1644.

\bibitem{MS02}
{\sc Myklebust, L.~I., and Skallerud, B.}
\newblock Model reduction methods for flexible structures.
\newblock In {\em Proc. 15th Nordic Seminar on Computational Mechanics,
  Aalborg, Denmark\/} (2002).

\bibitem{OAR89}
{\sc O'Callahan, J., Avitabile, P., and Riemer, R.}
\newblock System equivalent reduction expansion process {(SEREP)}.
\newblock In {\em Proc. 7. International Modal Analysis Conference, Las
  Vegas\/} (1989).

\bibitem{RS08}
{\sc Reis, T., and Stykel, T.}
\newblock Balanced truncation model reduction of second-order systems.
\newblock {\em Mathematical and Computer Modelling of Dynamical Systems 14}, 5
  (2008), 391--406.

\bibitem{SV07}
{\sc Salvini, P., and Vivio, F.}
\newblock Dynamic reduction strategies to extend modal analysis approach at
  higher frequencies.
\newblock {\em Finite Elements in Analysis and Design 43}, 11-12 (2007),
  931--940.

\bibitem{Schi2008}
{\sc Schilders, W. H.~A., van~der Vorst, H.~A., and Rommes, J.}, Eds.
\newblock {\em Model order reduction: Theory, research aspects and
  applications}, vol.~13 of {\em Mathematics in Industry}.
\newblock Spring{\-}er-Ver{\-}lag, Berlin, 2008.

\bibitem{SL07}
{\sc Tan, S. X.~D., and He, L.}
\newblock {\em Advanced Model Order Reduction Techniques for {VLSI} Designs}.
\newblock Cambridge University Press, 2007.

\bibitem{Ten04}
{\sc Teng, C.}
\newblock {\em Model Reduction of Second Order Linear Dynamical Systems}.
\newblock PhD thesis, Houston, Texas, 2004.

\bibitem{Van04}
{\sc Vandendorpe, A.}
\newblock {\em Model reduction of linear systems, an interpolation point of
  view}.
\newblock PhD thesis, Universit\'e catholique de Louvain, 2004.

\end{thebibliography}

\end{document}